\documentclass[12pt]{article}
\newtheorem{thm}{Theorem}
\newtheorem{cor}[thm]{Corollary}
\newtheorem{example}[thm]{Example}
\newtheorem{definition}[thm]{Definition}
\begin{document}

\title{Exceptional Prime Number Twins, Triplets and Multiplets}
\author{H. J. Weber\\Department of Physics\\
University of Virginia\\Charlottesville, 
VA 22904, U.S.A.}
\maketitle
\begin{abstract}
A classification of twin primes implies special twin 
primes. When applied to triplets, it yields exceptional 
prime number triplets. These generalize yielding 
exceptional prime number multiplets. 
\end{abstract}
\leftline{MSC: 11N05, 11N32, 11N80}
\leftline{Keywords: Twin primes, triplets, regular 
multiplets.} 


\section{Introduction}

Prime numbers have long fascinated people. Much is known 
about ordinary primes~\cite{hw}, much less about prime 
twins and mostly from sieve methods, and almost nothing 
about longer prime number multiplets. If primes and their 
multiplets are considered as the elements of the integers, 
then the aspects to be discussed in the following might 
be seen as part of a ``periodic table'' for them. 

When dealing with prime number twins, triplets and/or 
multiplets, it is standard practice to ignore as trivial 
the prime pairs $(2,p)$ of odd distance $p-2$ with $p$ 
any odd prime. In the following prime number multiplets 
will consist of odd primes only.   

\begin{definition}\label{1.1} A {\it generalized} twin 
prime consists of a pair of primes $p_i, p_f$ at distance 
$p_f-p_i=2D\geq 2.$ Ordinary twins are those at distance 
$2D=2.$ \end{definition}

\section{A classification of twin primes}

A triplet $p_i, p_m=p_i+2d_1, p_f=p_m+2d_2$ of prime 
numbers with $p_i<p_m<p_f$ is characterized by two 
distances $d_1, d_2.$ Each triplet consists of 
three generalized twin primes $(p_i, p_m), 
(p_m, p_f),\\(p_i, p_f)$. Empirical laws governing 
triplets therefore are intimately tied to those of 
the generalized twin primes. 

There is a basic parametrization of twin 
primes which generalizes to all prime number 
multiplets\cite{hjw1},\cite{hjw2}.  

\begin{thm}\label{2.1} Let $2D$ be the distance 
between two odd prime numbers $p_i,p_f$ of the pair, 
$D$ a natural number. There are three mutually 
exlusive classes of generalized twin primes. 
\begin{eqnarray}\nonumber 
I:~p_i&=&2a-D,~p_f=2a+D,~D~\rm{odd};\\\nonumber
II:~p_i&=&3(2a-1)-D,~p_f=3(2a-1)+D,~2|D,~3\not|D;
\\\nonumber 
III:~p_i&=&2a+1-D,~p_f=2a+1+D,~6|D,~D\geq 6,
\end{eqnarray}
where $a$ is the running integer variable. Values 
of $a$ for which a prime pair of distance $2D$ is 
reached are unpredictable (called arithmetic chaos). 

Only in class II are there special twins $3, 3+2D$ 
with a median $3+D$ not of the form $3(2a-1)$.  
\end{thm}
Each of these three classes contains infinitely 
many (possibly empty) subsets of prime pairs at 
various even distances.

{\it Proof.} Let us first consider the case of odd $D.$ 
Then $p_i=2a-D$ for some positive integer $a$ and, 
therefore, $p_f=p_i+2D=2a+D.$ The median $2a$ is the 
running integer variable of this class $I$. 

For even $D$ with $D$ not divisible by $3$ let 
$p_i=2n+1-D,$ so $p_f=2n+1+D$ for an appropriate integer 
$n.$ Let $p_i\neq 3,$ excluding a possible first pair 
with $p_i=3$ as {\it special.} Since one of three odd 
natural numbers at distance $D$ from each other is divisible 
by $3,$ the median $2n+1$ must be so, hence $2n+1=3(2a-1)$ 
for an appropriate integer $a\geq 2.$ Thus, the median 
$3(2a-1)$ of the pair $3(2a-1)\pm D$ is again a linear 
function of a running integer variable $a$. These prime 
number pairs constitute the 2nd class $II$. 

This argument is not valid for prime number pairs 
with $6|D,$ but these can obviously be parametrized as 
$2a+1\pm 6d,~D=6d.$ They comprise the 3rd and last class 
$III$ of generalized twins. Obviously these three classes 
are mutually exclusive and complete.$\diamond$ 
  
\begin{example}\label{1}  
Ordinary twins $2a\pm 1$ for $a=2, 3, 6, 9, 15, \ldots$ 
have $D=1$ and are in class $I$. For $D=3,$ the twins 
$2a\pm 3$ occur for $a=4, 5, 7, 8, 10, \ldots.$ For 
$D=5,$ the twins $2a\pm 5$ are at $a=4, 6, 9, 12, 18, 
\ldots.$   

For $D=2,$ the twins $3(2a-1)\pm 2$ for $a=2, 3, 4, 7, 
8, \ldots$ are in class $II.$ 

For $D=6,$ the twins $2a+1\pm 6$ for $a=5, 6, 8, 11, 12,
\ldots$ are in class $III.$ 
\end{example}

\begin{example}\label{3}
{\it Special twins} are the following prime number 
pairs $5\pm 2=(3,7); 7\pm 4=(3,11); 11\pm 8=(3,19);$ 
etc.\end{example}
How many special prime pairs are there? Since $3+2D$ 
forms an arithmetic progression with greatest common 
divisor $(3,2D)=1,$ Dirichlet's theorem~\cite{hw} 
on arithmetic progressions says there are infinitely 
many of them.

Theor.~\ref{2.1} tells us that all triplets are in 
precisely one of nine classes $(I,I),(I,II),(I,III),
(II,I),(II,II),(II,III),(III,I),\\(III,II),(III,III)$ 
which are labeled according to the class of their 
prime pairs. Quartets are distributed over $3^3$ 
classes, $n-$tuples over $3^{n-1}$ classes. 

It is well known that, except for the first pair $3, 5,$ 
ordinary twins all have the form $(6m-1, 6m+1)$ for some 
natural number $m.$ They belong to class $I.$ Using 
Theor.~\ref{2.1}, a second classification of generalized 
twins involving arithmetic progressions of conductor $6$ 
as their regular feature may be obtained~\cite{hjw1}. 

\begin{example}\label{2}.

Prime pairs at distance $2D=4$ are in class $II$ and of 
the form $6m+1, 6(m+1)-1$ for $m=1, 2, 3, \dots$ except 
for the singlet $3, 7.$ At distance $2D=6$ they are in 
class $I$ and have the form $6m-1, 6(m+1)-1$ for $m=1, 
2, 3, \ldots$ that are intertwined with $6m+1, 6(m+1)+1$ 
for $m=1, 2, 5, \dots.$ At distance $2D=12$ they are in 
class $III$ and of the form $6m-1, 6(m+2)-1$ for $m=1, 
3, 5, \ldots$ intertwined with $6m+1, 6(m+2)+1$ for $m=1, 
3, 5, \ldots.$ \end{example}
Clearly, the rules governing the form $6m\pm 1, 6m+b$ 
depend on the arithmetic of $D$ and $a.$

\section{Exceptional prime number triplets}

Let us start with some prime number triplet examples. 

\begin{example}\label{4} The triplet $~3, 5, 7$ is the 
only one at distances $[2, 2];$ the triplet $3, 7, 11$ 
is the only one at distances $[4, 4],$ and $3, 11, 19$ 
is the only one at distances $[8, 8],$ and $3, 13, 23$ 
is the only at $[10, 10],$ etc. Let's call such triplets 
at equal distances $[2D,2D]$ {\it exceptional or 
singlets.}\end{example} 

There are other classes of singlets. At distances 
$[2, 8],$ the prime triplet $3, 5, 13$ is 
the only one at distances $[2d_1=2,2d_2=8]$ and 
likewise is $3, 11, 13$ the only one at distances 
$[8, 2].$ 

{\it Rules for singlets or exceptions} among 
generalized triplet primes are the 
following\cite{hjw1}. 

\begin{thm}\label{3.1} (i) There is at most 
one generalized prime triplet with distances 
$[2D, 2D]$ for $D=1, 2, 4, 5, \ldots$ and 
$3\not| D.$  

(ii) When the distances are $[2d_1, 2d_2]$ with 
$3|d_2-d_1,$ and $3\not| d_1,$ there is at most 
one prime triplet $p_i=3, p_m=3+2d_1, p_f=
3+2d_1+2d_2$ for appropriate integers $d_1, d_2$.     
\end{thm}

{\it Proof.} (i) Because one of three odd numbers 
in a row at distances $[2D, 2D],$ with $D$ not 
divisible by $3,$ is divisible by $3,$ such a 
triplet must start with $3$. The argument fails 
when $3|D.$ 

(ii) Of $p_i,~p_m\equiv p_i+2d_1 \pmod{3},~
p_f\equiv p_i+4d_1 \pmod{3}$ at least one is 
divisible by $3,$ which must be $p_i.$~$\diamond$

Naturally, the question arises: Are there infinitely  
many such singlets, i.e. exceptional triplets? 
Dirichlet's theorem does not answer this one.  

\section{Exceptional prime number multiplets} 

It is obvious that {\it induced special multiplets} 
come about by adding a prime $p_f$ to a special twin 
$3,3+2D=p_m$ in class $II$ with $p_m-D=3+D
\neq 3(2a-1).$ This yields an {\it induced special 
triplet} $3,p_m=3+2D,p_f.$ Adding another prime 
generates an {\it induced special quartet}, etc. 

\begin{example}\label{5}
There are no prime multiplets at equal distance  
$2D,$\\$3 \not|D$ except triplets, because no other 
primes can be added at that distance, obviously. 
But there are quartets at equal distance $2D,3|D.$ 
They are in class $(I,I,I):~5,11,17,23;~41,47,53,
59;\\61,67,73,79;~251,257,263,269;~641,647,653,659;$ 
in class $III,\\III,III$ at equal distance $D$ with 
$6|D,$ etc. (Probably there is an infinite number 
of them.)\end{example} 

Thus, there are {\it no exceptional quartets at 
equal distance.} Likewise, there are {\it no 
exceptional $6-$tuples or $q-$tuples at equal 
distance, when $q$ has more than one prime 
divisor.}    

\begin{example}\label{9}
The {\it exceptional quintet} at equal distance 
is $5,11,17,\\23,29.$ There are no other such 
quintets, because one of $5$ odd numbers at 
distance $2D,3|D,5\not|D$ is divisible by $5.$ 

For $7$ the exceptional septet is $7,157,307,
457,607,757,907,$ at a surprisingly large 
distance. 

For the primes $11,13,19$ there are $6-$tuples 
$11,71,131,191,251,\\311;~11,491,971,1451,1931,
2411;~13,103,193,283,373,463;\\13,223,433,643,853,
1063;~13,3673,7333,10993,14653,18313;19,\\1669,
3319,4969,6619,8269$ and a $7-$tuple $17,2957,
5897,~8837,\\11777,14717,17657$ for $17,$ but 
the exceptional $11-,13-,17-, 19-$\\tuples are 
still at large, if they exist.\end{example}

These are special cases of the following  
rule for exceptional multiplets.  

\begin{cor}\label{3.3} {\it For any prime 
$p\geq 3$} there is at most an exceptional 
$p-$tuple $p,p+2D,\ldots, p+2(p-1)D$ at 
equal distance $2D,3|D,\\p \not|D.$\end{cor}  

{\it Proof.} A $p-$tuple composed of primes 
is exceptional, because one of any $p$ odd 
numbers at equal distance $2D,3|D,p\not|D$ is 
divisible by $p.~\diamond$  

When $q$ has several prime divisors, there are 
no new exceptional $q-$tuples at equal distance.   

Next we generalize item (ii) of Theor.~\ref{3.1} 
to quartets and subsequently higher multiplets. 

Exceptional quartets with $d_1\equiv d_2\pmod{3}, 
d_1\not\equiv 0\pmod{3},\\d_3\equiv 0\pmod{3}$ are 
possible that are {\it induced} by triplets, etc.

\begin{cor}\label{3.4} New exceptional quartets at 
distances $[2d_1,2d_2,2d_3]$\\when $d_2\equiv d_3
\pmod{3},d_2\not\equiv 0\pmod{3},d_1\equiv 0\pmod{3}$ 
and $d_1\equiv d_3\pmod{3},d_1\not\equiv 
0\pmod{3},d_2\equiv 0\pmod{3}.$\end{cor}  

{\it Proof.} This is clear from the proof of (ii) 
in Theor.~\ref{3.1}.~$\diamond$ 

The generalization to higher multiplets is 
straightforward and left to the reader. 

We conclude with an open question. Apart from 
exceptional multiplets, do all other multiplets 
repeat their pattern of given differences 
infinitely often?

\end{document}